\documentclass[12pt]{article}

\usepackage[T2A]{fontenc}
\usepackage[utf8]{inputenc}

\usepackage{amsfonts}
\usepackage{amssymb}
\usepackage{amsmath}
\usepackage{amsthm}

\def \le {\leqslant}
\def \ge {\geqslant}

\theoremstyle{plain}

\begin{document} 
 
\begin{Large}
\centerline{\bf Approximations to two real numbers}

\end{Large}
\vskip+0.5cm
\centerline{by  {\bf Igor D. Kan} and {\bf Nicky G. Moshchevitin}\footnote{research is supported by RFBF grant No. 09-01-00371a}
} 
\vskip2.0cm

\begin{small}

{\bf Abstract.}
Probably we have observed a
new simple phenomena dealing with approximations to two real numbers. 

\end{small}
\vskip+1.0cm
{\bf 1. The result.}

For a real $\xi$ denote the irrationality measure function
$$
\psi_\xi (t) = \min_{1\le x \le t}||x\xi ||.
$$
Here we suppose $x$ to be an integer number and $||\cdot ||$ stands for the distance to the nearest integer.

The main result of this note is the following

{\bf Theorem 1.}\,\,{\it
For any two different irrational 
numbers
$\alpha , \beta$
such that $\alpha\pm \beta \not\in \mathbb{Z}$ the difference function
$$
\psi_\alpha (t) - \psi_\beta (t)
$$
changes  its sign infinitely many times as $t\to +\infty$.}

The phenomenon observed in Theorem 1 cannot be generalized to
any dimension greater than one. In \cite{HIN}  the following two statements were proven.

{\bf Theorem 2.} (A. Khintchine, 1926) \,{\it
Let function $\psi (t)$ decreses to zero as $t\to +\infty$. Then there exist two algebraically independent
real numbers $\alpha^1,\alpha^2$ such that
for all $t$ large enough one has
$$
\psi_{(\alpha^1,\alpha^2)}
(t) :=
\min_{1\le \max(|x_1|,|x_2|) \le t} ||x_1\alpha^1+x_2\alpha^2||\le \psi (t). 
 $$}

{\bf Theorem 3.} (A. Khintchine, 1926) \,{\it
Let 
 $\psi_1(t)$ 
decreases to zero as  $t\to +\infty$ 
and  the function  $t\mapsto t\psi_1 (t)$ increases to infinity  as $t\to +\infty$.
Then there exist two algebraically independent
real numbers $\alpha_1,\alpha_2$ such that
for all $t$ large enough one has
$$
\psi_{
\left(
\begin{array}{c}
\alpha_1\cr \alpha_2
\end{array}\right)}
(t) : =
\min_{1\le x \le t} \, \max_{j=1,2}||x\alpha_j ||\le \psi_1 (t). 
 $$}

Of course in Theorems 2,3 we suppose $x_1,x_2,x$ to be integers.

Take  $\psi (t) = o(t^{-2}),\, t \to +\infty$.
 Take  $(\beta^1,\beta^2)$ to be numbers algebraically independent of $\alpha^1,\alpha^2$
such that they are badly approximable (in the sense of  a linear form):
$$
\inf_{(x_1,x_2) \in  \mathbb{Z}^2\setminus \{(0,0)\}} \left(
||x_1\beta^1+x_2\beta^2|| \cdot \max(|x_1|,|x_2|)^2 \right) >0.
$$
We see that
for all $t$ large enough  one has
$$
\psi_{(\alpha^1,\alpha^2)}
(t)<
\psi_{(\beta^1,\beta^2)}
(t).
$$

The similar situation holds in the case of simultaneous approximations.
Take  $\psi_1 (t) = o(t^{-1/2}),\, t \to +\infty$.
 Take  $\left(\begin{array}{c}\beta_1\cr \beta_2\end{array}\right)$ to be numbers algebraically independent of $\alpha^1,\alpha^2$
such that they are badly simultaneously approximable:
$$
\inf_{x \in  \mathbb{Z}\setminus \{0\}}
\,\left(\max_{j=1,2}
||x\beta_j|| \cdot |x|^{1/2} \right) >0.
$$
We see that

$$
\psi_{
\left(
\begin{array}{c}
\alpha_1\cr \alpha_2
\end{array}\right)}
(t) <
\psi_{
\left(
\begin{array}{c}
\beta_1\cr \beta_2
\end{array}\right)}
(t)
$$
for all $t$ large enough.
(Of course  here $ \psi_{(\beta^1,\beta^2)}, \psi_{
\left(
\begin{array}{c}
\beta_1\cr \beta_2
\end{array} \right) }$ are defined analogously to 
$ \psi_{(\alpha^1,\alpha^2)}, \psi_{
\left(
\begin{array}{c}
\alpha_1\cr \alpha_2
\end{array}\right)}$.)

{\bf 1. Proof of Theorem 1.}

We can assume that $0< \alpha,\beta <1$. 
We  consider continued fraction expansions
$$
\alpha = [0;a_1a_2,,...,a_n,...] , \,\,\, \beta= [0;b_1,b_2,...,b_n,...].
$$
Define
$$
\alpha_n = [a_n;a_{n+1}, a_{n+2},...], \,\,\,\, \alpha_n^* = [0; a_n,a_{n-1},...,a_1],
$$
$$
\beta_n = [b_n;b_{n+1}, b_{n+2},...], \,\,\,\, \beta_n^* = [0; b_n,b_{n-1},...,b_1],
$$
$$
\frac{r_n}{q_n} = [0;a_1,...,a_n],\,\,\,\,
\frac{s_n}{p_n} = [0;b_1,...,b_n].
$$

{\bf Lemma 1.}\,\,{\it For $n \ge 2$ one has
$$
||q_{n-1}\alpha ||q_{n+1} =
\frac{\alpha_{n+1} (a_{n+1}+\alpha_n^*)}{\alpha_{n+1}+\alpha_n^*}.
$$}

Proof. 

It is a well known fact
(see \cite{SCH}, Ch.1)
 that 
\begin{equation}\label{no}
\left|\alpha 
-\frac{r_{n-1}}{q_{n-1}}\right|=
\frac{1}{q_{n-1}^2(\alpha_n+\alpha_{n-1}^*)},
\end{equation}
 and
$$
\alpha_n^* = \frac{q_{n-1}}{q_n}.
$$
Instead of (\ref{no}) we can write
\begin{equation}\label{noo}
||q_{n-1}\alpha||=\frac{1}{q_{n-1}\alpha_n+q_{n-2}}.
\end{equation}

So we see that
$$
||q_{n-1}\alpha||q_{n+1} =
||q_{n-1}\alpha|| q_{n-1}
 \frac{q_{n}}{q_{n-1}}
 \frac{q_{n+1}}{q_n}
=
\frac{1}{(\alpha_n+\alpha_{n-1}^*)\alpha_n^*\alpha_{n+1}^*}
.
$$
But as
$$
\alpha_n = a_n +\frac{1}{\alpha_{n+1}},\,\,\,\,
a_n + \alpha_{n-1}^* = \frac{1}{\alpha_n^*}
$$
we see that
$$
\alpha_n+\alpha_{n-1}^* = \frac{1}{\alpha_n^*}+\frac{1}{\alpha_{n+1}}.
$$
So
$$
||q_{n-1}\alpha||q_{n+1} =
\frac{1}{\displaystyle{ \alpha_n^*\alpha_{n+1}^*\left(  
\frac{1}{\alpha_n^*}+\frac{1}{\alpha_{n+1}}
\right)}}=
\frac{\alpha_{n+1}}{\alpha_{n+1}^*(\alpha_n^*+\alpha_{n+1})} =
\frac{\alpha_{n+1} (a_{n+1}+\alpha_n^*)}{\alpha_{n+1}+\alpha_n^*}
.
$$
Lemma is proved.

As $a_{n+1} \ge 1$ and $ \alpha_{n+1}>1$ we obtain the following

{\bf Corollary.}\,\,{\it 
For $n \ge 2$ one has
\begin{equation}\label{0}
||q_{n-1}\alpha||q_{n+1} >1.
\end{equation}}

{\bf Lemma 2.}\,\,{\it
Suppose that $m,n\ge 2$ and 
\begin{equation}\label{1}
q_{n+1} \le p_{m+1} .
\end{equation}
Then 
\begin{equation}\label{2}
||q_{n-1}\alpha ||> ||p_m\beta ||.
\end{equation}}

Proof.

Suppose that (\ref{2}) is not true. Then from (\ref{1}) and (\ref{0}) we see that
$$
1< ||q_{n-1}\alpha||q_{n+1}\le 
||p_{m}\beta||p_{m+1}
.$$
As 
(see \cite{SCH}, Ch.1)
$$
||p_{m}\beta ||p_{m+1} 
=
\frac{1}{\displaystyle{1 +\frac{\beta_{m+1}^*}{\beta_{m+2}}}}<1
$$
we have a contradiction.
Lemma 2 is proved.

Now we are able to prove theorem 1.

Consider the sequences    
$$
q_0\le q_1<...<q_n<q_{n+1}<...,\,\,\,\, p_0\le p_1<...,<p_m<p_{m+1}<...
$$
of convergents' denominators to  $\alpha, \beta$
correspondingly.
Suppose that the statement of theorem 1 is false for certain   irrationalities $\alpha,\beta$. Without loss of generality assume that
for all $t \ge p_{m_0}\ge q_{n_0-1}$ one has
\begin{equation}\label{con}
\psi_\beta (t) \ge \psi_\alpha (t).
\end{equation}
From Lemma 2 and the asumption (\ref{con}) we see that between two consecutive denominators $p_m,p_{m+1}, m \ge m_0$ 
not more than one denominator of the form $q_n$ may occur. Here we give a proof of this fact.
Let $q_{n-1}\le p_m<q_n<q_{n+1}<...<q_{n+t}\le p_{m+1}$ and $t\ge 1$.
Then
$$
||p_m\beta ||= \psi_\beta (p_m) \ge \psi_\alpha (p_m) =\psi_\alpha (q_{n-1}) =
||q_{n-1}\alpha||
$$
and
$$
q_{n+1}\le q_{n+t}<p_{n+1}.
$$
This contradicts to Lemma 2.

So we can define the  sequence of integers
$$
m_0\ge 1,\,\,\, m_j \ge m_{j-1}+1
$$
such that
$$
p_{m_0}<q_{n_0}\le p_{m_0+1}<...<p_{m_1}<q_{n_0+ 1}\le
p_{m_1+1}< ...<p_{m_2}< q_{n_0+1}\le p_{m_2+1}< ...
$$
$$
 < p_{m_{j-1}}< q_{n_0+j-1} \le p_{m_{j-1}+1} <
... < p_{m_{j}}< q_{n_0+j} \le p_{m_j+1} <...<p_{m_{j+1}}< q_{n_0+j+1} \le p_{m_{j+1}+1}< ... .
$$
By (\ref{con})
we see that  for all $j\ge 0$ 
  one has
  \begin{equation}\label{two}
||q_{n_0+j-1}\alpha ||=\psi_\alpha (q_{n_0+j-1})
=\psi_\alpha (p_{m_j}) \le\psi_\beta (p_{m_j}) =||p_{m_j}\beta ||.
\end{equation}
From (\ref{con}) we also have
\begin{equation}\label{one}
||q_{n_0+j}\alpha||=\psi_\alpha (q_{n_0+j})
=\psi_\alpha (p_{m_j+1}) \le\psi_\beta (p_{m_j+1}) =||p_{m_j+1}\beta||.
\end{equation}

We distinguish two cases.
In the {\bf first case} we suppose that for infinitely many $j$    at least one of  the inequalities in 
(\ref{two},\ref{one}) is strict, that is there is the sign $<$ insead of $\le$. In the {\bf second case} for all $j$
large enough  we have equalities in both     (\ref{two},\ref{one}).

Consider the {\bf first case}.  Without loss of generality we assume that
\begin{equation}\label{two1}
||q_{n_0+j-1}\alpha||=\psi_\alpha (q_{n_0+j-1})
=\psi_\alpha (p_{m_j}) <\psi_\beta (p_{m_j}) =||p_{m_j}\beta||.
\end{equation}
 From (\ref{noo}) we have
$$
 ||q_{n_0+j-1}\alpha||=
\frac{1}{q_{n_0+j-1}\alpha_{n_0+j}+q_{n_0+j-2}},
\,\,\,\,
||p_{m_j}\beta||=
\frac{1}{p_{m_j}\beta_{m_j+1}+p_{m_j-1}}.
$$
So
\begin{equation}\label{part}
p_{m_j}\beta_{m_j+1}+p_{m_j-1} <
q_{n_0+j-1}\alpha_{n_0+j}+q_{n_0+j-2}
\end{equation}
As
$$
\beta_{m_j+1}=b_{m_j+1}+\frac{1}{\beta_{m_j+2}},\,\,\,\,\,
\alpha_{n_0+1}=a_{n_0+j}+\frac{1}{\alpha_{n_0+j+1}}
$$
from  (\ref{part}) we deduce that
$$
p_{m_j}\left( b_{m_j+1}+\frac{1}{\beta_{m_j+2}}\right)<
q_{n_0+j-1}\left(a_{n_0+j}+\frac{1}{\alpha_{n_0+j+1}}\right) +q_{n_0+j-2}
$$
or
$$
p_{m_j+1}+\frac{p_{m_j}}{\beta_{m_j+2}} <
q_{n_0+j}+\frac{q_{n_0+j-1}}{\alpha_{n_0+j+1}}.
$$
But
$$
{p_{m_j+1}}\ge q_{n_0+j},\,\,\,\, p_{m_j} \ge q_{n_0+j-1}.
$$
So
\begin{equation}\label{oro}
\beta_{m_j+2} > \alpha_{n_0+j-1}.
\end{equation}

From the other hand from (\ref{one}) we deduce that
$$
\frac{1}{q_{n_0+j}\alpha_{n_0+j+1} +q_{n_0+j-1} }=
||q_{n_0+j}\alpha|| =\psi_\alpha (q_{n_0+j}) =
\psi_\alpha (p_{m_j+1}) \le
$$
$$
\le 
\psi_\beta (p_{m_j+1}) =
||p_{m_j+1}\beta||=
\frac{1}{p_{m_j+1}\beta_{m_j+2}+p_{m_j}}
.
$$
So
$$
p_{m_j+1}\beta_{m_j+2}+p_{m_j}\le q_{n_0+j}\alpha_{n_0+j+1} +q_{n_0+j-1} .
$$
As
$$
q_{n_0+j-1} \le p_{m_j},\,\,\,\,\,
q_{n_0+j} \le p_{m_j+1} 
$$
we see that
$$
 \beta_{m_j+2}\le \alpha_{n_0+j+1}.
$$
This contradicts (\ref{oro}).

In the {\bf second case}  we see that for $j$
 large enough  one has
$$\psi_\beta (p_{m_j+1} ) =\psi_\alpha (q_{{n_0}+j})=
\psi_\beta (p_{m_{j+1}}).
$$
Hence
$$
m_j+1 = m_{j+1}.
$$
 But in the case under consideration 
 we see that there exist  $m_0,n_0$  such that
$$
p_{m_0+j}\beta -q_{n_0+j}\alpha = \pm r_{n_0+j}\pm s_{m_0+j},
$$
 $$
p_{m_0+j+1}\beta -q_{n_0+j+1}\alpha = \mp r_{n_0+j+1}\mp s_{m_0+j+1},
$$
 where the choise of the signs $\pm$ depends on the lenghts of the corresponding continued fractons.
Remind that $\alpha,\beta $ are irrational numbers. So
$$
p_{m_0+j}q_{n_0+j+1}-
p_{m_0+j+1}q_{n_0+j}=0
$$
and
 $$
\frac{p_{m_0+j}}{p_{m_0+j+1}} =
[0;b_{m_0+j+1},b_{m_0+j},...,b_1] =
\frac{q_{n_0+j}}{q_{n_0+j+1}} =
[0;a_{m_0+j+1},a_{m_0+j},...,a_1] , \,\,\,\, j =1,2,3,...
$$
and so $\alpha =\pm\beta$.

The proof of Theorem 1 is complete.

\end{document}